\documentclass[10pt,twoside]{siamart1116}

\usepackage[english]{babel}
\usepackage{graphicx,epstopdf,epsfig}
\usepackage{amsfonts,epsfig,fancyhdr,graphics, hyperref,amsmath,amssymb}
\usepackage{slashbox}
\graphicspath{ {./Images/} }
\newtheorem{remark}[theorem]{Remark}
\newtheorem{example}[theorem]{Example}

\renewcommand{\theequation}{\thesection.\arabic{equation}}
\renewtheorem{theorem}{Theorem}[section]

\newcommand{\Title}{A Note on Generalized Locally Toeplitz Operators}
\begin{document}
	\setcounter{page}{1}
	
	\thispagestyle{empty}
	\title{\Title}
	\author{
		V. B. Kiran Kumar\thanks{Department of Mathematics, Cochin University of Science And Technology, Kalamassery, Kochi, Kerala 682022
			(kiranbalu36@gmail.com)}
		\and
		N. S. Sarathkumar\thanks{Department of Mathematics, Cochin University of Science And Technology, Kalamassery, Kochi, Kerala 682022 ( sarathkumarns05@gmail.com).}}
	\maketitle
	\begin{abstract}
		Generalized Locally Toeplitz (GLT)  matrix sequences arise from large linear systems that approximate Partial Differential Equations (PDEs),
		Fractional Differential Equations (FDEs), and Integro-Differential Equations (IDEs).  GLT sequences of matrices have been developed
		to study the spectral/singular value behaviour of the numerical approximations to various PDEs, Fades and IDEs.  These approximations can be achieved
		using any discretization method on appropriate grids through local techniques such as Finite Differences, Finite Elements,
		Finite Volumes, Isogeometric Analysis, and Discontinuous Galerkin methods. Spectral and singular value symbols are essential for analyzing the eigenvalue and singular value distributions of
		matrix sequences in the Weyl sense. In this article, we provide a comprehensive overview of the operator-theoretic aspect of GLT sequences. The theory of GLT sequences, along with findings on the asymptotic spectral distribution of perturbed matrix sequences, is a highly effective and successful method for calculating the spectral symbol $f$. Therefore, developing an automatic procedure to compute the spectral symbols of these matrix sequences would be advantageous, a task that Ahmed Ratnani, N S Sarathkumar, S. Serra-Capizzano have partially undertaken. As an application of the theory developed here, we propose an automatic procedure for computing the symbol of the underlying sequences of matrices, assuming they form a GLT sequence that meets mild conditions.
	\end{abstract}
	\section{Introduction}
	In order to model practical issues, Differential Equations (DEs) are frequently used in physics, engineering, and applied sciences. The analytical solution of such DEs typically lacks a closed form. Therefore, it is crucial to approximate the  solution $u$ of DE using a numerical method.
	
	If we take into account a linear DE
	 $$Au = f $$and a linear numerical approach, the actual calculation of the numerical solution reduces to solving a linear system $$A_nu_n = f_n$$ whose size $d_n$ diverges with $n$.
	 
	To compute accurate solutions in an acceptable amount of time, it is essential to solve high-dimensional linear systems efficiently. In this regard, it is well known that the convergence properties of popular iterative solvers, such as multigrid and preconditioned Krylov techniques, are heavily dependent on the spectral features of the matrices to which they are applied. In order to select or create the approximate solution and discretization method, it is helpful to be aware of the asymptotic spectral distribution of the sequence $\{A_n\}_n$. This assertion is supported by the fact that significant estimates of the superlinear convergence of the conjugate gradient method produced by Beckermann and Kuijlaars in \cite{beckermann2001superlinear} are closely related to the asymptotic spectral distribution of the considered matrices.
	
	When the eigenvalues of such matrices represent physically relevant parameters, the spectral distribution of DE discretization matrices is crucial itself in addition to the construction and analysis of suitable solvers.
	 This is true for a large group of issues that come up in engineering and applied sciences, like the research of natural vibration frequencies for an elastic material; for more information, see review \cite{garoni2019symbol} and the references therein.
	In addition, it is frequently noticed in practise that the spectral distribution of the sequence $\{A_n\}_n$ is somehow linked to the spectrum of the differential operator $A$.
	
	The theory of GLT sequences is an effective tool for computing/analyzing the asymptotic spectral distribution of the discretization matrices resulting from the numerical approximation of continuous problems, such as  Partial Differential Equations (PDEs), Fractional Differential Equations (FDEs), Integro Differential Equations (IDEs). When discretizing a
	linear PDE by means of a linear numerical method, the actual computation of the numerical solution reduces to solving a linear system
	$A_nu_n = b_n.$ The size $d_n$ of this linear system increases when the discretization parameter $n$ tends to infinity. Hence, what we actually
	have is not just a single linear system, but a whole sequence of linear systems with increasing size. For example, if we consider the Schr\"odinger operator that maps $f\mapsto -f^{''}+vf$, where $v$ is a real-valued periodic potential function, then the corresponding finite difference approximation leads to a sequence of block Toeplitz matrices. If we consider more general differential operators like those arising from diffusion problems ($f\mapsto (-af^{'})^{'}+vf$) or convection-diffusion-reaction problems  ($f\mapsto (-af^{'})^{'}+bf^{'}+vf$), we end up with Locally Toeplitz (LT) or Generalized Locally Toeplitz (GLT) sequences \cite{garoni2017generalized}. What is often observed in
	practice is that the sequence of discretization matrices $\{A_n\}_n$ enjoys an asymptotic spectral distribution, which is somehow related to
	the spectrum of the differential operator associated with the considered PDE. 
	
	Tilli's research on locally Toeplitz (LT) sequences \cite{Tilli-LT} and the spectral theory of Toeplitz matrices serve as the foundation for the theory of GLT sequences. These days, the books \cite{garoni2017generalized,garoni2018generalized,barbarino2020block,barbarinomulti} are the primary and finest sources in the literature for theory and applications of GLT spaces. In these volumes, we may get an extensive and comprehensive description of GLT sequences, block GLT sequences, and their respective multivariate counterparts. Numerous studies investigated possible generalisations of such findings due to the inherent relationship between matrix sequences and the related symbol functions. They attempted to identify symbol functions or function spaces corresponding to specific classes of matrix sequences in order to comprehend the spectral asymptotic. Researchers like S. Serra-Capizzano, A. Bottcher, G. Barbarino, C. Garoni, etc. have recently started such studies in the context of GLT sequences \cite{barbarino2017equivalence, barbarino2017convergence,barbarino2020block,Bottcher,garoni2017generalized,sarath}.
	
	These advancements mostly concentrated on matrix symbols and sequences. The operator theoretic version of GLT sequences is almost untouched. The primary motive of this article is to provide an overview of the operator related to GLT sequences and investigate its characteristics. In section 2 we give some preliminaries for the article. We construct a Hilbert space in section 3 to define the GLT operator.  Additionally, we introduce two fundamental operators on this Hilbert space that are analogous  Toeplitz and multiplication operators, which play a key role in the definition of the GLT operator. Along with the operator norm, we also provide two convergence notions, and we examine which convergence notion is most suitable for defining the GLT operator. LT operator and GLT operator are introduced in section 4, using that suitable convergence notion. Also we  investigate some theories of GLT operators that analogies to GLT sequences. The theory of GLT sequences, coupled with the results on the asymptotic spectral distribution of perturbed matrix sequences, is among the most powerful and successful tools for computing the spectral symbol $f$. Consequently, developing an automatic procedure to compute the spectral symbols of these matrix sequences would be advantageous, a goal that Ahmed Ratnani has partially pursued. n \cite{Sarath2024}, an automatic procedure was proposed for computing the symbol of the underlying sequences of matrices, assuming it is a GLT sequence that satisfies mild conditions. In Section 5, we establish a technique, derived from operator theory, to compute the spectral symbol of GLT sequences, accompanied by numerical illustrations.
	
	\section{Preliminaries}
	In practice, it is often observed that the sequence of discretization matrices $\{A_n\}_n$ exhibits an asymptotic spectral distribution. The concepts of asymptotic spectral distribution and asymptotic singular value distribution are closely related.
	\begin{definition}
		\normalfont
		Let $\{A_n\}_n$ be a matrix-sequence with size n for each $A_n$ and $f:D\subset \mathbb{R}^k \to \mathbb{C}$ be a measurable function.
		\begin{itemize}
			\item We say that $\{A_n\}_n$ has an asymptotic singular value distribution with symbol $f$, and  write
			$\{A_n\}_n\sim_{\sigma} f$, if, for all $F\in C_c(\mathbb{R})$,  the space of complex-valued continuous functions defined on $\mathbb{R}$ and with compact support,
			$$ \lim _{n\to\infty} \frac{1}{n}\sum_{i=1}^{n}F(\sigma_i(A_n)) = \frac{1}{\mu_k(D)}\int_D F(|f(x)|) dx, $$
			where $\sigma_i(A_n),\, i=1,\cdots,n,$ are the singular values of $A_n$, $\mu_k$ is the Lebesgue measure in $\mathbb{R}^k.$
			\item We say that $\{A_n\}_n$ has an asymptotic  eigenvalue distribution with symbol $f$, and  write
			$\{A_n\}_n\sim_\lambda f$, if, for all $F\in C_c(\mathbb{C})$,  the space of complex-valued continuous functions defined on $\mathbb{C}$ and with compact support,
			$$ \lim _{n\to\infty} \frac{1}{n}\sum_{i=1}^{n}F(\lambda_i(A_n)) = \frac{1}{\mu_k(D)}\int_D F(f(x))dx, $$
			where $\lambda_i(A_n),\, i=1,\cdots,n,$ are the eigenvalues of $A_n$.
		\end{itemize}
		
	\end{definition} 
	The concept of approximating classes of sequences, introduced in \cite{Serraacs}, forms the basis of a spectral approximation theory for matrix sequences. It provides tools for computing the asymptotic spectral or singular value distribution of a matrix sequence $\{A_n\}_n$.
	\begin{definition}\normalfont
		Let $\lbrace A_{n}\rbrace _{n} $ be a matrix-sequence and $\lbrace\lbrace B_{n,m} \rbrace_n\rbrace_m$ a sequence of matrix-sequences. We say that $\lbrace\lbrace B_{n,m} \rbrace_n\rbrace_m$ is an approximating class of sequences $\mathrm{(a.c.s)}$ for  $\lbrace A_{n}\rbrace _{n}$ if the following condition is met: for every $m$ there exists an $n_m$ such that, for $n \geqslant n_m$,
		
		$$A_n= B_{n,m}+R_{n,m}+N_{n,m}$$
		
		$\mathrm{rank} (R_{n,m}) \leqslant c(m)n,\quad \| N_{n,m} \| \leqslant \omega(m) $ 
		where $\parallel.\parallel$ is the spectral norm,  $ n_m, c(m)$ and $\omega(m)$ depend only on $m$ and 
		\begin{center}
			$\lim_{m\to\infty} c(m) = \lim_{m\to\infty} \omega(m)=0$.
		\end{center}
	\end{definition}

	Now, following \cite{garoni2017generalized}, we give the construction of the GLT matrix sequences, originally defined in \cite{capizzano2003generalized,serra2006glt}.
	\begin{definition}
		\normalfont Let $ a:[0,1]\rightarrow \mathbb{C}$ be a Riemann-integrable function and $f \in L^1{([-\pi,\pi])}$. We say that a matrix-sequence $\lbrace A_{n}\rbrace _{n}$ is a \textit{Locally Toeplitz (LT) sequence} with symbol $a \otimes f$, and we write $\lbrace A_{n}\rbrace _{n}\sim_{\mathrm{LT}} a\otimes f$, if $\lbrace \lbrace LT_{n}^{m}(a,f) \rbrace_n\rbrace_{m \in \mathbb{N}}$ is an a.c.s. for $\lbrace A_{n}\rbrace _{n}$, where
		\begin{equation*}
			\begin{aligned}
				LT_{n}^{m}(a,f)&=[D_m(a) \otimes T_{\lfloor n/m \rfloor}(f)]\oplus O_{n (mod\, m)}\\
				&=\underset{i=1,...,n}{\mathrm{diag}}[a\left(\frac{i}{m}\right)T_{\lfloor n/m \rfloor }(f)]\oplus O_{n (mod\, m)}.
			\end{aligned}
		\end{equation*}
		Here $T_n(f)$ is the Toeplitz matrix generated by the function $f$, $O_{n (mod\, m)}$ is the zero matrix of order $n (mod\, m)$,  and $D_m(a)$  is a  $(m \times m)$ diagonal matrix associated with $a$ given by
		$$D_m(a)=\underset{i=1,...,n}{\mathrm{diag}} a\left(\frac{i}{m}\right).$$
	\end{definition}
	For example, consider the second order differential problem $-(a(x)u^\prime(x))^\prime=f(x),\quad x\in (0,1),$ with initial condition 
	$u(0)=\alpha,\quad u(1)=\beta$, where $a:[0,1]\to \mathbb{R}$ is continuous. Then by descretizing this problem using classical second-order central FD scheme, we will get a linear system $A_nu_n=\left(\frac{1}{n+1}\right)^2f_n$, where
	\begin{equation*}
		A_n=\begin{Bmatrix}
			a(x_{\frac{1}{2}})+a(x_{\frac{3}{2}})&-a(x_{\frac{3}{2}})&&&&\\
			-a(x_{\frac{3}{2}})&a(x_{\frac{3}{2}})+a(x_{\frac{5}{2}})&-a(x_{\frac{5}{2}})&&&\\
			&-a(x_{\frac{5}{2}})&a(x_{\frac{5}{2}})+a(x_{\frac{7}{2}})&-a(x_{\frac{7}{2}})&&\\
			&&-a(x_{\frac{7}{2}})&\ddots&\ddots&\\
			&&&\ddots&\ddots&-a(x_{n-\frac{1}{2}})\\
			&&&&-a(x_{n-\frac{1}{2}})&a(x_{n-\frac{1}{2}})+a(x_{n+{\frac{1}{2}}})
		\end{Bmatrix}.
	\end{equation*}
Now the matrix sequence $\{A_n\}_n$ is a LT sequence with symbol $a\otimes f$.
	
	\begin{definition}\label{GLT defn}\normalfont
		Let  $\kappa:[0,1] \times [-\pi,\pi] \rightarrow \mathbb{C}$ be a measurable function. We say that a matrix-sequence $\lbrace A_{n}\rbrace _{n}$  is a \textit{Generalized  Locally Toeplitz (GLT) sequence} with symbol $\kappa$, and we write $\lbrace A_{n}\rbrace _{n}\sim_{\mathrm{GLT}} \kappa$, if the following conditions is hold.\\
		For every $m$ varying in some infinite subset of $\mathbb{N} $ there exists a  finite number of LT sequences $\lbrace A_{n}^{(i,m)}\rbrace _{n} \sim_{\mathrm{LT}} a_{i,m}\otimes f_{i,m},i=1,...,k_m ,$ such that:
		\begin{itemize}
			\item $\sum \limits_{i=1}^{k_m} a_{i,m} \otimes f_{i,m}\rightarrow \kappa $ in measure over $[0,1] \times [-\pi,\pi]$ when $m \rightarrow  \infty$;
			\item $\left\lbrace\left\lbrace \sum \limits_{i=1}^{k_m} A_{n}^{(i,m)}\right\rbrace _{n}\right\rbrace_m $ is an a.c.s. for $\lbrace A_{n}\rbrace _{n}$.
		\end{itemize}
	\end{definition}

As a result of the findings in \cite{barbarino2018normal}, the subsequent Lemma offers a proof for the relationship between the operators defined in the following sections and GLT sequences.
\begin{lemma}\label{nGLT}\cite{barbarino2018normal}
	Let $a\in C^\infty([0,1])$ and $f$ be a trigonometric polynomial, then $LT_n^{\lfloor n\rfloor}(a,f)$ is a Locally Toeplitz sequence with symbol $a\otimes f$.
\end{lemma}
A seminorm $q$ on the space of all matrix sequences is introduced in \cite{sarath}  motivated by the notion of a.c.s. This seminorm plays a crucial role in defining the GLT operator.
$$q(\{A_\textbf{\textit{n}}\}_n):=\inf \left\{\limsup_{n\to \infty}\frac{\|N_\textbf{\textit{n}}\|_{2}}{N(\textbf{\textit{n}})^{1/2}}:\,\, R_\textbf{\textit{n}}+N_\textbf{\textit{n}}=A_\textbf{\textit{n}}, \,\, {\text{rank} (R_\textbf{\textit{n}})}=o(N(\textbf{\textit{n}}))\right\},$$
Here the infimum is taken over all such decompositions of $A_\textbf{\textit{n}}$ and $\|N_\textbf{\textit{n}}\|_2$ is the Frobenius  norm of $N_\textbf{\textit{n}}$.
	A class of matrix-sequences that plays a central role in the framework of the theory of GLT sequences is the class of zero-distributed sequences.
\begin{definition}\label{zero}\normalfont\cite{sarath}
	A matrix-sequence $\{Z_n\}_n$ is said to be a \textit{zero-distributed sequence} if $q(\{Z_n\}_n)=0.$
\end{definition}
The following lemma is a consequence of the results in \cite{tyrtyshnikov1996unifying} and \cite{sarath},
\begin{lemma}\label{zerod}
	Let $\{Z_n\}_n$ be a matrix sequence. If $\|Z_n\|_{2}^2=o(n)$, then $\{Z_n\}$ is zero distributed sequence.
\end{lemma}
 We recall some basic properties  of  GLT matrix-sequences from \cite{barbarino2017equivalence,garoni2017generalized};
\begin{theorem}\label{tog1}
	Let $\{A_n\}_n\sim_{GLT}f$ and $\{B_n\}_n\sim_{GLT}g$. Then,
	\begin{enumerate}
		\item $\{A_{n}^{*}\}_n\sim_{GLT}\bar{f}$.
		\item$\{\alpha A_n+\beta B_n\}_n\sim_{GLT}\alpha f+\beta g$, for all $\alpha, \beta \in \mathbb{C}.$
		\item $\{A_nB_n\}_n\sim_{GLT} fg$.
		\item if $\{A_n\}_n\sim_{GLT} h$ then $f=h$ a.e.
		\item $\{A_n\}_n$ is zero-distributed iff $f=0$ a.e.
	\end{enumerate}
\end{theorem}
The notion of approximating class of sequences have a fundamental role in the theory of GLT sequences. In this section, we introduce a convergence notion analogous to a.c.s for defining GLT operator.
Let $I=(0,1]\cap \mathbb{Q}$. 
$$\mathcal{L}^2= \bigoplus_{l\in I}L_j=\{(f_l)\in \prod_{l\in I}L_l : \sum_{l\in I}\|f_l\|^{2} < \infty\}$$
where $L_l=  L^2[-\pi,\pi]$
$$ \mathcal{H}^2= \bigoplus_{l\in I}H_i=\{(f_l)\in \prod_{l\in I} H_l : \sum_{l\in I}\|f_l\|^{2} < \infty\}$$
where $H_l=  H^2[-\pi,\pi]$.\\
We define $e_{lk}:I\to \mathbb{C}$ as $e_{l,k}(j)=\delta_{lj}e^{ikt}$. Then $\mathcal{B_L}=\{e_{lk}: k\in \mathbb{Z} , l\in I\}$ and $\mathcal{B_H}=\{e_{lk}: k\in \mathbb{Z}^+ , l\in I\}$ will form orthonormal basis for $\mathcal{L}^2$ and $\mathcal{H}^2$ respectively. Note that $\mathcal{L}^2$ and $\mathcal{H}^2$ form a Hilbert space with innerproduct $\sum_{l\in I}<f_l,g_l>$. Let $L(\mathcal{H}^2)$ be the space of all linear operator whose domain contains $B_\mathcal{H}$. In next section we will introduce some topologies.	
	\section{Topologies on some subspaces of $L(\mathcal{H}^2)$}

	An essentially bounded measurable function $f$ on $[-\pi,\pi]$ and a Riemann Integrable function $a$ on $[0,1]$ induces two operators,  on $\mathcal{L}^2$ and  on $\mathcal{H}^2$ respectively as follows. The weighted Multiplication Operator $M_{af}$ is
	$$ M_{af}((g_l))=(a(l)fg_l)$$
	for every $(g_l)$ in $\mathcal{L}^2$. The (Locally) Toeplitz Operator $T_{af}$ is defines in terms of the orthogonal projection $P$ from $\mathcal{L}^2$ onto $\mathcal{H}^2$, as
	$$T_{af}((g_l))=PM_{af}$$
	for every $(g_l)$ in $\mathcal{H}^2$.\\
	
	Some standard properties of Multiplication  and Toeplitz operators are also enjoyed by $M_{af}$ and $T_{af}$, since the proofs are near imitation (see \cite{martinez2007introduction} for example), we omit them.\\
	If $a: [0,1]\to \mathbb{C}$ be a Riemann Integrable function and  $ f\in L^\infty$, then
	\begin{enumerate}
		\item  $\|M_{af}\|=\|a\|_{\infty}\|f\|_{\infty}$.
		\item $\sigma(M_{af})=ess\,ran (a\otimes f)$.
		\item The Spectral Inclusion Theorem : The spectrum of $M_{af}$ is contained in the spectrum of $T_{af}$. More precisely, $ess\, ran(a\otimes f)=\Pi(M_{af})=\sigma(M_{af})\subset \Pi(T_{af})\subset \sigma(T_{af})$, where $\Pi$ is the approximate spectrum.
		\item $r(T_{af})=\|T_{af}\|=\|M_{af}\|=\|a\|_\infty\|f\|_\infty$.
	\end{enumerate}
	Now we can see the connection between $T_{af}$ and LT sequences. We see that wew can construct LT sequences from $t_{af}$. For $m,n\in \mathbb{N} $ and $n\geq m$, $P_{n}^{m}$ be the projection operator on $\mathcal{H}^2$ to the subsapce spanned by $\{e_{lk}:k=0,1,2,\ldots,\lfloor\frac{n}{m}\rfloor-1,l\in I_m\}$, where $I_m=\{\frac{j}{m} : j=1,2,\cdots,m\}$. Here onwards we use the notation  $P_n$ instead of $P_n^{\lfloor n\rfloor}$.
	
	\begin{theorem}
		Let  $a: [0,1]\to \mathbb{C}$ be a Riemann Integrable  function  and $ f\in L^\infty[-\pi,\pi]$ with Fourier series $\sum_{n=-\infty}^{\infty}b_ne^{in\theta} $, then the matrix representation  $T_n^m$ of the operator $P_n^mT_{af}P_n^m$ is given by 
		$$T_n^m=D_m(a) \otimes T_{\lfloor n/m \rfloor}(f)=\underset{i=1,...,m}{\mathrm{diag}}a(\frac{i}{m})T_{\lfloor n/m \rfloor (f)}.$$
	\end{theorem}
	\begin{proof}
		$T_n^m$ is the matrix of order $m\cdot\lfloor\frac{n}{m}\rfloor$ and its $(j,k)^{th}$ entry $(T_n^m)_{jk}$ is given by $ \langle T_{af}e_{j_1j_2},e_{k_1k_2}  \rangle$, where $j_1=\frac{1}{m}\frac{j}{\lfloor \frac{n}{m}\rfloor}, j_2=j (mod\, \lfloor \frac{n}{m}\rfloor), k_1=\frac{1}{m}\frac{k}{\lfloor \frac{n}{m}\rfloor}, k_2=k (mod\, \lfloor \frac{n}{m}\rfloor).$\\
		Now 	$\langle T_{af}e_{j_1j_2},e_{k_1k_2}\rangle= 0$ for $l_1\neq k_1$ and 
		$$\langle T_{af}e_{j_1k_1},e_{j_1k_2}\rangle=\frac{1}{2\pi}\int_{0}^{2\pi}a(j_1)fe^{i(j_2-k_2)t} =a(j_1)f_{k_2-j_2}.$$
		.	
	\end{proof}	
	Now we  define two bounded linear operator $M_{jk}:\mathcal{L}^2 \to \mathcal{L}^2$ and $T_{jk}:\mathcal{H}^2\to \mathcal{H}^2$, which are the special cases of $M_{af}$ and $T_{af}$. \\
	Let $M_{jk},j,k\in\mathbb{Z}$ be the weighted multiplication operator on $\mathcal{L}^2$ defined by 
	$$ M_{jk}((f_l))=(e^{\pi i(2jl+kt)}f_l)$$
	Also $T_{jk}$ be a bounded linear operator on $ \mathcal{H}^2$ defined by
	$$ T_{jk}=PM_{jk}P$$
	For $n\in \mathbb{N}$, let
	
	$T_{n_{jk}}=P_nT_{jk}P_n +O_{n(mod\lfloor n\rfloor)}$. Then $\{T_{n_{jk}}\}_n$ is a GLT sequence with symbol $e^{2\pi i jx}\cdot e^{ikt}$. Let ${E}$ be the space spaned by the linear operators $T_{jk};j,k\in \mathbb{Z}$. These $T_{jk}$s are the building blocks for GLT operators.\\
	
	The notion of approximating class of sequences plays a key role in the definition of Locally Toeplitz sequences as well as in the GLT sequences. So we need to introduce a convergence notion analogous to the approximating class of sequences to define LT operator and GLT Operator.
	
	\subsection{Operator norm}
	Let $\bar{{E}} $ is the closure of ${E}$ in $B(\mathcal{H}^2)$ with respect to the operator norm.
	\begin{theorem}
		Let $T\in \bar{E}$ and $T_n=P_nTP_n$. Then $\{T_n\}_n$ is GLT sequence.
	\end{theorem}
	\begin{proof}
		Let $T\in \bar{E}$, then there exists a sequence of operators $\{T_m\}_m$ in $E$ such that $\lim_{m\to \infty}\|T-T_m\|=0$, where $T_m=\sum_{j=j_{m}}^{j_{m\prime}}\sum_{k=k_{m}}^{k_{m\prime}}a_{jk}T_{jk}$. Now,
		$$\lim_{m\to \infty}\|P_n(T-T_m)P_n\|=0.$$
		That is $\{P_nT_mP_n\}_n$ is an a.c.s. for $\{A_n\}_n$ and by Lemma \ref{nGLT} ${\{P_nT_mP_n\}_n}\sim_{LT}\sum_{j=j_{m}}^{j_{m\prime}}\sum_{k=k_{m}}^{k_{m\prime}}a_{jk}e^{2\pi i jx}\cdot e^{ikt}$. Then $\{A_n\}_n$ is a GLT sequence.
	\end{proof}	
\textbf{Advantages:}
\begin{itemize}
	\item The truncations $\{T_n\}_n$ of the operator $T\in \bar{E}$ is GLT sequence.
	\item Let $T,S\in \bar{E}$ and $\{T_n\}\sim_{GLT} f$,$\{S_n\}_n\sim_{GLT} g$, then $\alpha T+\beta S\in \bar{E}$, since $\bar{E}$ is vector space. By Theorem \ref{tog1} $\{\alpha T_n+\beta S_n\}_n \sim_{GLT} {\alpha f+\beta g}$.
	\item Let $T\in\bar{{E}}$ and $\{T_n\}_n\sim_{GLT} f$, then $T^*\in \bar {E}$ and by Theorem \ref{tog1} $\{T^*_n\}_n\sim_{GLT} \bar{f}$.
\end{itemize}
\textbf{Disadvantages:}
	\begin{itemize}
		\item The symbol of GLT sequences not necessarily bounded. For example  consider $\{T_n(f)\}_n$ , a sequence of Toeplitz matrices with an unbounded symbol $f$. Such GLT sequences do not arise from an element in $\bar{{E}}$, since it consist  of only bounded operators.
		\item From the theory of GLT sequences, we expect that the truncation of all compact operators is  GLT. But the case may be different if we consider the operator norm convergence. Since only zero operator give rise to a zero distributed operator.  
	\end{itemize}	
	Motivated by the notion of seminorm $q$ in \cite{sarath} , we introduce a norm on ${E}$ and a semi innerproduct  on a subspace of $L(\mathcal{H}^2)$.

\subsection{An innerporoduct on ${E}$}	
Consider $T_1,T_2\in {E}$. We define an  innerproduct
 $\langle\cdot,\cdot\rangle_G:{E}\times {E}\to \mathbb{C}$, such that
\begin{equation}\label{innerprrdct}
	\langle T_1,T_2 \rangle_G=\lim_{n\to\infty}\frac{1}{n}\sum_{r=0}^{\lfloor \frac{n}{\lfloor \sqrt{n}\rfloor}\rfloor-1}\sum_{l\in I_{\lfloor \sqrt{n}\rfloor}}\langle T_1(e_{lr}),T_2(e_{lr})\rangle.
\end{equation}
Also $\|\cdot\|_G$ be the norm induced by the innerproduct $\langle\cdot,\cdot\rangle_G$. Let $E^{\prime}$ be the completion of ${E}$ with respect to the norm induced from the innerproduct.
\begin{theorem}
	$\{T_{jk}:j,k\in \mathbb{Z}\}$ is an orthonormal basis of $E$.
\end{theorem}
\begin{proof}	
	\begin{equation*}
		\begin{aligned}
			\langle T_{j_1k_1},T_{j_2k_2}\rangle_G&=\lim_{n\to\infty}\frac{1}{n}\sum_{r=0}^{\lfloor \frac{n}{\lfloor \sqrt{n}\rfloor}\rfloor-1}\sum_{l\in I_{\lfloor \sqrt{n}\rfloor}}\frac{1}{2\pi}\int_{-\pi}^{\pi}T_{j_1k_1}(e_{lr})\cdot \overline{T_{j_2k_2}(e_{lr}) }dt,\\
			&=\lim_{n\to\infty}\frac{1}{n}\sum_{r=0}^{\lfloor \frac{n}{\lfloor \sqrt{n}\rfloor}\rfloor-1}\sum_{l\in I_{\lfloor \sqrt{n}\rfloor}}\frac{1}{2\pi}\int_{-\pi}^{\pi}(e^{2\pi ij_1l}e^{i(k_1+r)t}) \overline{(e^{2\pi ij_2l}e^{i(k_2+r)t})}dt,\\
			&=\lim_{n\to\infty}\frac{1}{n}\sum_{r=0}^{\lfloor \frac{n}{\lfloor \sqrt{n}\rfloor}\rfloor-1}\sum_{l\in I_{\lfloor \sqrt{n}\rfloor}}\frac{1}{2\pi}\int_{-\pi}^{\pi}(e^{2\pi i(j_1-j_2)l}e^{i(k_1-k_2)t})dt	\\
			&=\left\{	\begin{aligned}
				0\quad\quad\quad\quad\quad\quad\quad\quad	&\text{if $k_1\neq k_2,$}\\
				\lim_{n\to\infty}\frac{1}{n}\sum_{r=0}^{\lfloor \frac{n}{\lfloor \sqrt{n}\rfloor}\rfloor-1}\sum_{l\in I_{\lfloor \sqrt{n}\rfloor}}e^{2\pi i(j_1-j_2)l}\,\,\,\,&\text{if $k_1 = k_2$.}	
			\end{aligned}\right.\\
			\langle T_{j_1k_1},T_{j_2k_2}\rangle_G&=\lim_{n\to\infty}\frac{1}{n}\sum_{r=0}^{\lfloor \frac{n}{\lfloor \sqrt{n}\rfloor}\rfloor-1}\sum_{l\in I_{\lfloor \sqrt{n}\rfloor}}e^{2\pi i(j_1-j_2)l},
				\end{aligned}
		\end{equation*}
	\begin{equation*}
\begin{aligned}
			&=\lim_{n\to\infty}\frac{1}{\lfloor \sqrt{n}\rfloor}\sum_{l\in I_{\lfloor \sqrt{n}\rfloor}}e^{2\pi i(j_1-j_2)l},\\
			&=\int_{0}^{1}e^{2\pi i(j_1-j_2)lx}dx,\\
			&=\left\{	\begin{aligned}
				0\quad\quad	&\text{if $j_1\neq j_2,$}\\
				1\quad\quad&\text{if $j_1 = j_2$.}
			\end{aligned}\right.
		\end{aligned}
	\end{equation*}
\end{proof}
 The following theorem says that the truncation of all the elements in $E^{\prime}$ is a GLT sequences.
\begin{theorem}\label{GLT}
	If $T\in E^{\prime}$, then $\{P_nTP_n\oplus O_{n \mod \lfloor \sqrt{n}\rfloor}\}_n$ is a GLT sequence.
\end{theorem}
\begin{proof}
Let ${T}_n=\{P_nTP_n\oplus O_{n \mod \lfloor \sqrt{n}\rfloor}\}_n$ and ${T}_n^m=\{P_nT_mP_n\oplus O_{n \mod \lfloor \sqrt{n}\rfloor}\}_n$.	For $\epsilon>0$, there exist an integer $M$ such that, for all $m>M$,
	$$\|T- \sum_{j=j_m,k=k_m}^{j_m^\prime,k_m^\prime}a_{jk}T_{jk}\|_G<\epsilon.$$
	\begin{equation*}
		\begin{aligned}
			q_{2}^{2}(\{{T}_n-{T}_{n,m}\}_n)&\leq \lim\frac{1}{n}\|{T}_n-{T}_{n,m}\|_F^2\\
			&\leq \lim_{n\to\infty}\frac{1}{n}\sum_{r=0}^{\lfloor \frac{n}{\lfloor \sqrt{n}\rfloor}\rfloor-1}\sum_{l\in I_{\lfloor \sqrt{n}\rfloor}}\|T- \sum_{j,k}a_{jk}T_{jk}\|^2\\
			&<\epsilon.
		\end{aligned}
	\end{equation*}	
Now $\{{T}_n^m\}_n$ is an a.c.s. for $\{{T}_n\}_n$ and by Theorem \ref{nGLT} ${\{{T}_n^m\}_n}\sim_{LT}\sum_{j=j_m,k=k_m}^{j_m^\prime,k_m^\prime}a_{jk}e^{2\pi i jx}\cdot e^{ikt}$. Then $\{{T}_n\}_n$ is a GLT sequence.
\end{proof}
From Theorem \ref{GLT} we understood that, if $T$ is in $E^{\prime}$ such that $\sum_{j_m,k_m}a_{j_mk_m}T_{j_mk_m}$ converges to $T$, then there exist a measurable function $\kappa$ such that $\sum_{j_m,k_m}a_{j_mk_m}e^{2\pi j_mx}\cdot e^{ik_mt}$ converges to $\kappa$ in $L^p$ and $\{T_n\}_n$ is a GLT sequence with symbol $k$. \\
The next theorem says that the truncation of Principal Locally Toeplitz operator is a GLT sequences.
\begin{theorem}\label{LT-GLT}
	$LT_{af}$ belongs to $E^{\prime}$. In other words sequence $\{{T}_n\}_n$ is a GLT sequence.
\end{theorem}
\begin{proof}
	Let $\sum_{j,k}a_{jk}e^{2\pi ijx}e^{ikt}$	converges to $a\otimes f$ in $L^2$ norm. Let $T=T_{af}-\sum_{j,k}a_{jk}T_{jk}$. Now consider
\begin{equation*}
		\begin{aligned}
			\|T\|_G^2&=\lim_{n\to\infty}\frac{1}{n}\sum_{r=0}^{\lfloor \frac{n}{\lfloor \sqrt{n}\rfloor}\rfloor-1}\sum_{l\in I_{\lfloor \sqrt{n}\rfloor}}\frac{1}{2\pi}\int_{-\pi}^{\pi}|a(l)\cdot e^{irt}\cdot f-\sum_{j,k}a_{jk}e^{2\pi ijl}\cdot e^{i(k+r)t}|^2 dt,\\
		&=\lim_{n\to\infty}\frac{1}{n}\sum_{r=0}^{\lfloor \frac{n}{\lfloor \sqrt{n}\rfloor}\rfloor-1}\sum_{l\in I_{\lfloor \sqrt{n}\rfloor}}\frac{1}{2\pi}\int_{-\pi}^{\pi}\left(|a(l)f|^2-\sum_{j,k}a(l)f\cdot \overline{(a_{jk})}e^{-i\pi (2jl+kt)}\right.\\
		&\left.\quad-\overline{a(l)f}\sum_{j,k}a_{jk}e^{i\pi(2jl+kt)}-\left|\sum_{j,k}a_{jk}e^{2\pi ijl}\cdot e^{i(k+r)t}\right|^2\right).
		\end{aligned}
\end{equation*}
\begin{equation*}
		\begin{aligned}
			\lim_{n\to\infty}\frac{1}{n}\sum_{r=0}^{\lfloor \frac{n}{\lfloor \sqrt{n}\rfloor}\rfloor-1}\sum_{l\in I_{\lfloor \sqrt{n}\rfloor}}\frac{1}{2\pi}\int_{-\pi}^{pi}\|\sum_{j,k}a_{jk}e^{2\pi ijl}\cdot e^{i(k+r)t}|^2|a(l)f|^2&=\|f\|_2^2\limsup_{n\to \infty}\frac{1}{\sqrt{n}}\sum_{l\in I_{\lfloor \sqrt{n}\rfloor}}|a(l)|^2,\\
			&=\|f\|_2^2\|a\|_2^2.
		\end{aligned}
\end{equation*}
\begin{equation*}
		\begin{aligned}
			\lim_{n\to\infty}\frac{1}{n}\sum_{r=0}^{\lfloor \frac{n}{\lfloor \sqrt{n}\rfloor}\rfloor-1}\sum_{l\in I_{\lfloor \sqrt{n}\rfloor}}\frac{1}{2\pi}\int_{-\pi}^{\pi}&\sum_{j,k}a(l)f\cdot \overline{(a_{jk})}e^{-i\pi (2jl+kt)}\\
			&=\sum_{j,k}\lim_{n\to\infty}\frac{1}{\lfloor \sqrt{n}\rfloor}\sum_{l\in I_{\lfloor \sqrt{n}\rfloor}}\frac{a(l)}{2\pi}\int_{-\pi}^{\pi}a(l)f\cdot \overline{(a_{jk})}e^{-i\pi (2jl+kt)},\\
			&=\sum_{j,k}\lim_{n\to\infty}\frac{1}{\lfloor \sqrt{n}\rfloor}\sum_{l\in I_{\lfloor \sqrt{n}\rfloor}}\frac{a(l)}{2\pi}\int_{-\pi}^{\pi}a(l)(\sum_{m}a_me^{imt})\cdot \overline{(a_{jk})}e^{-i\pi (2jl+kt)},\\
			&=\sum_{j,k}\lim_{n\to\infty}\frac{1}{\lfloor \sqrt{n}\rfloor}\sum_{l\in I_{\lfloor \sqrt{n}\rfloor}}a(l)e^{-2\pi i jl}a_k\overline{a_{jk}},\\
			&=\sum_{j,k}a_k\overline{_{jk}}\int_{0}^{1}a(x)e^{-2\pi i jx}dx,
		\end{aligned}
	\end{equation*}
\begin{equation*}
\begin{aligned}
			&=\sum_{j,k}a_k\overline{_{jk}}\int_{0}^{1}(\sum_{p}a_pe^{2\pi ipx})e^{-2\pi i jx}dx,\\
			&=\sum_{j,k}\overline{a_jk}a_ka_j,\\
			&=\sum_{j,k}|a_{jk}|^2.
		\end{aligned}
\end{equation*}
	Similarly,
	$$\lim_{n\to\infty}\frac{1}{n}\sum_{r=0}^{\lfloor \frac{n}{\lfloor \sqrt{n}\rfloor}\rfloor-1}\sum_{l\in I_{\lfloor \sqrt{n}\rfloor}}\frac{1}{2\pi}\int_{-\pi}^{\pi}\sum_{j,k}\overline{a(l)f}\cdot {(a_{jk})}e^{-i\pi (2jl+kt)}=\sum_{j,k}|a_{jk}|^2.$$
	$$\lim_{n\to\infty}\frac{1}{n}\sum_{r=0}^{\lfloor \frac{n}{\lfloor \sqrt{n}\rfloor}\rfloor-1}\sum_{l\in I_{\lfloor \sqrt{n}\rfloor}}\frac{1}{2\pi}\int_{-\pi}^{\pi}|\sum_{j,k}a_{jk}e^{2\pi ijl}\cdot e^{i(k+r)t}|^2=\sum_{j,k}|a_{jk}|^2.$$
	Then,
$$	\|T\|_G^2= \|f\|_2^2\|a\|_2^2-\sum_{j,k}|a_{jk}|^2<\epsilon.$$
\end{proof}

\textbf{Advantages:}
\begin{itemize}
	\item The truncations $\{{T}_n\}_n$ of the operator $T\in\bar{E}$ is GLT sequence.
	\item Let $T,S\in E^{\prime}$ and $\{T_n\}\sim_{GLT} f$,$\{S_n\}_n\sim_{GLT} g$, then $\alpha T+\beta S\in E^{\prime}$ and $\{\alpha T_n+\beta S_n\}_n \sim_{GLT} {\alpha f+\beta g}$.
	\item $\bar{E}\subset E^{\prime}$,and $E^{\prime}$ is much richer than $E^{\prime}$.
\end{itemize}
\textbf{Disadvantages:}
	\begin{itemize}
		\item Product of two operators in ${E}^{\prime}$ may not be in ${E}^{\prime}$. The example is given in Remark \ref{example}.
		\item If $T\in E^{\prime}$, then $T^*$ may not belongs to $E^{\prime}$. For example,
		\begin{equation*}
			T(e_{lr})=\left\{
			\begin{aligned}
				e_{lr}\quad&\text{if } r=0,\\
				0\quad\quad&\text{otherwise}.
			\end{aligned}\right.
		\end{equation*}
	Then $\|T\|_G=0$ and $\|T\|_G=1$, contradicting the fact that if $T$ and $T^*$ belongs to $E^\prime$, then $\|T\|_G=\|T^*\|_G$. 
	\end{itemize}
\subsection{A semi norm on a subspace of  $L(\mathcal{H}^2)$ }
	\begin{definition}
		Let $T\in L(\mathcal{H}^2)$ and define a function $Q:L(\mathcal{H}^2)\to \mathbb{R}$ as
		$$Q(T)=\inf\left\{\limsup_{n\to \infty}\frac{1}{\sqrt{n}}\left(\sum_{e_{lk}\in C_n}\|T(e_{lk})\|^2\right)^{\frac{1}{2}}: C_n\sqcup D_n=B_n, \lim_{n\to\infty}\#(D_n)=0\right\},$$
		where $B_n=\{e_{lk}:k=0,1,2,\ldots,\lfloor\frac{n}{\sqrt{n}}\rfloor-1,l\in I_{\lfloor\sqrt{n}\rfloor}\}$, $\#(S)$ is the cardinality of $S$ and infimum is taken over all such decomposition of $B_n$.
	\end{definition}
	Define the subspace ${L}_{g}$ of $L(\mathcal{H}^2)$  as follows 
	$$ 
	L_{g}=\left\{T\in L(\mathcal{H}^2)\,\,:\,\,Q(T)<\infty\right\}.$$

	Let $L_G=L_g/Z(\mathcal{H}^2)$ be the quotient space of $L_g$, where  $Z(\mathcal{H}^2)$ be the space of all operators $Z$ with $Q(Z)=0.$
	\begin{theorem}
		$L_G$ is a Hilbert space.
	\end{theorem} 
	\begin{proof}
		Let $T_1,T_2 \in L_g$. Then for every $m\in \mathbb{N}$, there exist four set of vectors $C_n^1, D_n^1,C_n^2,D_n^2$ such that $C_n^1\sqcup D_n^1=C_n^2\sqcup D_n^2=B_n,$ and
		$$\limsup_{n\to\infty}\frac{1}{\sqrt{n}}\left(\sum_{e_{lk}\in C_n^i}\|T_i(e_{lk}\|^2)\right)^{\frac{1}{2}}\leq Q(T_i)+\frac{1}{m},\quad \forall i=1,2.$$
		Now we verify the axioms of seminorm. The non negativity and is trivial. For triangular inequality,
		Let us take $C_n=C_n^1\cap C_n^2$ and $D_n=D_n^1\cup D_n^2$. Then
		$$\lim_{n\to\infty}\frac{1}{n}\#(D_n)\leq\lim_{n\to\infty}\frac{1}{n}\#(D_n^1)+\lim_{n\to\infty}\frac{1}{n}\#(D_n^2)=0.$$
		Now,
		\begin{equation*}
			\begin{aligned}
				Q(T_1+T_2)&\leq \limsup_{n\to\infty}\frac{1}{\sqrt{n}}\left(\sum_{e_{lk}\in C_n}\|(T_1+T_2)(e_{lk}\|^2)\right)^{\frac{1}{2}}\\
				&\leq \limsup_{n\to\infty}\frac{1}{\sqrt{n}}\left(\sum_{e_{lk}\in C_n^1}\|T_1(e_{lk}\|^2)\right)^{\frac{1}{2}}+\limsup_{n\to\infty}\frac{1}{\sqrt{n}}\left(\sum_{e_{lk}\in C_n^2}\|T_2(e_{lk}\|^2)\right)^{\frac{1}{2}}\\
				&\leq Q(T_1)+Q(T_2)+\frac{2}{m}.
			\end{aligned}
		\end{equation*} 
		Thus $\|T_1+T_2\|_g\leq \|T_1\|_g+\|T_2\|g.$\\
		$\|\alpha T\|_g=|\alpha|\|T\|_g$ is straight forward.
		Hence $Q$ is a seminorm in $L_g$. Then the function $\tilde{Q}:L_G\to \mathbb{R}$ defined as 
		$$\tilde{Q}(T+Z)=Q(T),$$
		becomes a norm on $L_G$.\\
		Next we prove the completeness of $L_G$. Let $\{T_m\}$ be a Cauchy sequence in $L_G$.  It suffice to show  that the convergence of a subsequence. Then there exist a subsequence say $\{T_m\}$ itself such that 
		$$\tilde{Q}(T_{k+1}-T_k+Z)< \frac{1}{2^k}\quad\forall k.$$
		Fix a number $k$. Then there exist decomposition $\{C_n^k\}$, $\{D_n^k\}$ of $\{B_n\}$ and increasing sequence of positive integers $\{N_k\}$, such that 
		$$\frac{1}{n}\left(\sum_{e_{lk}\in C_n^k}\|(T_{k+1}-T_k)(e_{lk}\|^2)\right)< \frac{1}{4^k} \quad \forall n>N_k, \text{ and }
		\lim_{n\to\infty}\frac{1}{n}\#(D_n^\prime)=0,$$
		where $\displaystyle D_n^\prime =\bigcup_{i=0}^jC_n^{k+i}$ if $N_{k+j}\leq N_{k+j+1}$.\\
		Let $E_n^k=\{e_{lr}\in C_n^k: \|(T_{k+1}-T_k)(e_{lr})\|\geq\frac{1}{\sqrt{2}^k}\}.$ Then
		\begin{equation*}
			\begin{aligned}
				\frac{1}{2^k}\limsup_{n\to\infty}\frac{1}{n}\#(E_n^k)&\leq\limsup_{n\to\infty}\frac{1}{n}\sum_{e_{lr}\in E_n^k}\|(T_{k+1}-T_k)(e_{lr})\|^2\\
				&\leq \limsup_{n\to\infty}\frac{1}{n}\sum_{e_{lr}\in C_n^k}\|(T_{k+1}-T_k)(e_{lr})\|^2\\
				&\leq \frac{1}{2^{2k}}.
			\end{aligned}
		\end{equation*}
		Then,
		$$\limsup_{n\to\infty}\frac{1}{n}\#(E_n^k)\leq\frac{1}{2^{k}}\quad\forall k.$$
		Now consider the sets 
		$$F^\prime=\{e_{lr}\in \mathcal{B}_{\mathcal{H}}: \|(T_{k+1}-T_k)(e_{lr})\|\geq \frac{1}{k/2} \text{ for infinitly many k}\}, \text{ and}$$
		$$F_n^\prime =\{e_{lr}\in F^\prime:r=0,1,2,\ldots,\lfloor\frac{n}{\sqrt{n}}\rfloor-1,l\in I_{\lfloor\sqrt{n}\rfloor}\}.$$
		Then,
		$$\limsup_{n\to\infty}\frac{1}{n}\#(F_n^\prime)=\limsup_{n\to\infty}\frac{1}{n}\#(F_n^\prime\cup {D_n^k})\leq \limsup_{n\to\infty}\frac{1}{n}\#(E_n^k\cup D_n^k)=\limsup_{n\to\infty}\frac{1}{n}\#(E_n^k)\leq\frac{1}{2^k}.$$
		The above inequality is true for every $k$, therefore $\lim_{n\to\infty}\frac{1}{n}\#(|F_n^\prime)=0.$
		Let $F=\mathcal{B}_{\mathcal{H}}\backslash F^\prime$.
		Then for every  $e_{lr}$ in $ F$, there exist an index $k_{lr}$ such that $$\|(T_{k+1}-T_k)(e_{lr})\|<\frac{1}{2^{k/2}}\quad \forall k\geq k_{lr}.$$
		
		Now define an operator $T$ as
		\begin{equation*}
			T(e_{lr})=\left\{
			\begin{aligned}
				\lim_{k\to\infty}T_k(e_{lr})\quad&\text{if } e_{l,r}\in F,\\
				0\quad\quad\quad&\text{otherwise}.
			\end{aligned}\right.
		\end{equation*}
		Let,
		$$G_{k+j}=\left\{e_{lr}\in  F:\|(T_{k+j+i+1}-T_{k+j+i})(e_{lr})\|<\frac{1}{2^{\frac{k+j+i}{2}}}\,\forall i \text{ and }\|(T_{k+j}-T_{k+j-1})(e_{l,r})\|>\frac{1}{2^\frac{k+j-1}{2}}\right\}.$$
		Then $G_{k+j_1}\subseteq G_{k+j_2}$, whenever $j_1\leq j_2$.\\
		Consider $\displaystyle C_n=F_n\bigcap_{i=0}^{j}C_n^{k+i}, \text{ if } N_{k+j}\leq n <N_{k+j+1}.$ Then $\#(B_n\backslash C_n)=o(n)$. Now,
		\begin{equation*}
			\begin{aligned}
				{\tilde{Q}}^2(T-T_k)&\leq \limsup_{n\to\infty}\frac{1}{n}\sum_{e_{lr}\in C_n}\|(T-T_k)(e_{lr})\|^2\\
				&\leq\limsup_{n\to\infty}\frac{1}{n}\sum_{e_{lr}\in C_n\cap G_{k+j_n}}\|(T-T_k)(e_{lr})\|^2+\sum_{e_{lr}\in C_n\backslash G_{k+j_n}}\|(T-T_k)(e_{lr})\|^2,
			\end{aligned}
		\end{equation*}
		\begin{equation*}
			\begin{aligned}
		\limsup_{n\to\infty}\frac{1}{n}\sum_{e_{lr}\in C_n\cap G_{k+j_n}}\|(T-T_k)(e_{lr})\|^2		&= \limsup_{n\to\infty}\frac{1}{n}\sum_{e_{lr}\in C_n\cap G_{k+j_n}}\|(T-T_{k+j_n}+T_{k+j_n}-T_k)(e_{lr})\|^2,\\
				&\leq\limsup_{n\to \infty}\frac{1}{n}\sum_{e_{lr}\in C_n\cap G_{k+j_n}}\left(\|
				(T-T_{k+j_n})(e_{lr})\|^2+ \|(T_{k+j_n}-T_k)(e_{lr})\|^2\right.\\
				&+ \left.2\|(T-T_{k+j_n})(e_{lr})\|\|(T_{k+j_n}-T_k)(e_{lr})\|\right)\\
				&\leq \limsup\frac{1}{n}\left(\sum_{i=1}^\infty2^{\frac{-1}{k+j_n+i}}+\sum_{i=1}^\infty 2^{\frac{-1}{k+i}}+2\sum_{i=1}^\infty2^{-(\frac{k+j_n+i}{2})}\sum_{i=1}^\infty 2^{-\frac{k+i}{2}}\right)\\
				&=0.
			\end{aligned}
		\end{equation*}
		For the Hilbert space structure, let $\tilde{Q}(T_1+Z)=\alpha$ and $\tilde{Q}(T_2+Z)=\beta$. Then foe every $\epsilon >0$, there exist sequences of sets $\{C_{1,n}\},\{C_{2,n}\},\{D_{1,n}\},\{D_{2,n}\}$ and a positive integer $n_0$ such that $\#(D_{1,n})=\#(D_{2,n})=o(n)$, and for all $n>n_0$, 
		$$ \frac{1}{n}\sum_{e_{lr}\in C_{1,n}}\|T_1(e{lr})\|^2 \leq \alpha ^2+\epsilon,\quad\quad
		\frac{1}{n}\sum_{e_{lr}\in C_{2,n}}\|T_2(e_{lr})\|^2 \leq \beta ^2+\epsilon.$$
		Let $C_n=C_{1,n}\cap C_{2,n}$.	Now,
		$$Q^2(T_1+T_2)\leq \limsup_{n\to\infty}\frac{1}{n}\sum_{e_{lr}\in C_n}\|(T_1+T_2)(e_{lr})\|^2,\quad
		Q^2(T_1-T_2)\leq \limsup_{n\to\infty}\frac{1}{n}\sum_{e_{lr}\in C_n}\|(T_1-T_2)(e_{lr})\|^2.$$
		Then,
		\begin{equation*}
			\begin{aligned}
				Q^2(T_1+T_2)+ Q^2(T_1-T_2)&\leq \limsup_{n\to\infty}\frac{1}{n}\sum_{e_{lr}\in C_n}\left(\|(T_1+T_2)(e_{lr})\|^2+ \|(T_1-T_2)(e_{l,r})\|^2\right)\\
				&\leq 2\alpha^2 +2 \beta^2+ 4\epsilon.
			\end{aligned}
		\end{equation*}
		Thus,
		\begin{equation}\label{parallel1}
			Q^2(T_1+T_2)+ Q^2(T_1-T_2)\leq 2Q^2(T_1)+2Q^2(T_2).
		\end{equation}
		Let $Q(T_1+T_2)=\alpha^\prime$ and $Q(T_1-T_2)=\beta^\prime$. Then for every $\epsilon$, there exist sequences of sets $\{C_{1,n}^\prime\}$,$\{C_{2,n}^\prime\}$,$\{D_{1,n}^\prime\}$,\\$\{D_{2,n}^\prime\}$ and a positive integer $n_1$ such that $\#(D_{1,n}^\prime)=\#(|D_{2,n}^\prime)=o(n)$ and for all $n>n_1$, 
		$$ \frac{1}{n}\sum_{e_{lr}\in C_{1,n}^\prime}\|(T_1+T_2)(e_{lr})\|^2 \leq {\alpha^\prime} ^2+\epsilon,\quad\quad
		\frac{1}{n}\sum_{e_{lr}\in C_{2,n}^\prime}\|(T_1-T_2)(e_{l,r})\|^2 \leq {\beta^\prime} ^2+\epsilon.$$
		Let $C_n^\prime=C_{1,n}^\prime\cap C_{2,n}^\prime$. Now
		\begin{equation*}
			\begin{aligned}
				4Q^2(T_1)&\leq \limsup_{n\to\infty}\sum_{e_{lr}\in C_n^\prime}\|(T_1+T_2+T_1-T_2)(e_{l,r})\|^2,\\
				4Q^2(T_2)&\leq \limsup_{n\to\infty}\sum_{e_{lr}\in C_n^\prime}\|(T_1+T_2-T_1+T_2)(e_{l,r})\|^2,\\
				4Q^2(T_1)+4Q^2(T_2)&\leq \limsup_{n\to\infty}\frac{1}{n}\sum_{e_{lr}\in C_n^\prime}\left(\|(T_1+T_2+T_1-T_2)(e_{lr})\|^2+ \|(T_1+T_2-T_1+T_2)(e_{lr})\|^2\right)\\
				&\leq \limsup_{n\to\infty}\frac{2}{n}\left(\sum_{C_n^\prime}\|(T_1+T_2)(e_{lr})\|^2+\sum_{C_n^\prime}\|(T_1-T_2)(e_{lr})\|^2 \right)\\
				&\leq 2{\alpha^\prime}^2+2{\beta^\prime}^2+4\epsilon.
			\end{aligned}
		\end{equation*}
		Thus,\begin{equation} \label{parallel2}
			Q^2(T_1+T_2)+ Q^2(T_1-T_2)\geq 2Q^2(T_1)+2Q^2(T_2).
		\end{equation}
	From (\ref{parallel1}) and (\ref{parallel2})
		
	$$Q^2(T_1+T_2)+ Q^2(T_1-T_2)= 2Q^2(T_1)+2Q^2(T_2).$$
	Since the norm $\tilde{Q}$ satisfies parallelogram equality, $L_G$ forms a Hilbert space. Let  $\mathcal{G}$ be the closure of $E$ in $L_G$ with respect to the metric induced from $Q$.
	
	\end{proof}
	
\begin{theorem}\label{onb}
	$\{T_{jk}+Z:j,k\in \mathbb{Z}\}$ is an orthonormal basis of $\mathcal{G}$.
\end{theorem}
\begin{proof}
	\begin{equation*}
		\begin{aligned}
			\tilde{Q}^2(T_{j_1k_1}+T_{j_2k_2}+Z)&={Q}^2(T_{j_1k_1}+T_{j_2k_2})\\
			&=\limsup_{n\to \infty}\frac{1}{n}\sum_{e_{l,r}\in C_n}\|(T_{j_1k_1}+T_{j_2k_2})(e_{l,r})\|^2\\
			&=\limsup_{n\to \infty}\frac{1}{n}\sum_{r=0}^{\lfloor \frac{n}{\lfloor \sqrt{n}\rfloor}\rfloor-1}\sum_{l\in I_{\lfloor \sqrt{n}\rfloor}}\sum_{l\in I}\frac{1}{2\pi}\int_{-\pi}^{\pi}|(T_{j_1k_1}+T_{j_2k_2})(e_{l,r})|^2\\
				&=\limsup_{n\to \infty}\frac{1}{n}\sum_{r=0}^{\lfloor \frac{n}{\lfloor \sqrt{n}\rfloor}\rfloor-1}\sum_{l\in I_{\lfloor \sqrt{n}\rfloor}}\frac{1}{2\pi}\int_{-\pi}^{\pi}|e^{2\pi ij_1}\cdot e^{i(k_1+r)t}+e^{2\pi ij_2}\cdot e^{i(k_2+r)t}|^2\\
			&=\left\{	\begin{aligned}
				2\quad\quad\quad\quad\quad\quad\quad\quad	&\text{if $k_1\neq k_2,$}\\
				\lim_{n\to\infty}\frac{1}{n}\sum_{r=0}^{\lfloor \frac{n}{\lfloor \sqrt{n}\rfloor}\rfloor-1}\sum_{l\in I_{\lfloor \sqrt{n}\rfloor}}|e^{2\pi ij_1l}+e^{2\pi i j_2l}|^2\,\,\,\,&\text{if $k_1 = k_2$.}	
			\end{aligned}\right.
			\end{aligned}
	\end{equation*}
\begin{equation*}
\begin{aligned}
			{Q}^2(T_{j_1k}+T_{j_2k})&=	\lim_{n\to\infty}\frac{1}{n}\sum_{r=0}^{\lfloor \frac{n}{\lfloor \sqrt{n}\rfloor}\rfloor-1}\sum_{l\in I_{\lfloor \sqrt{n}\rfloor}}|e^{2\pi ij_1l}+e^{2\pi i j_2l}|^2\\
			&=\lim_{n\to\infty}\frac{1}{\lfloor \sqrt{n}\rfloor}\sum_{l\in I_{\lfloor \sqrt{n}\rfloor}}|e^{2\pi ij_1l}+e^{2\pi i j_2l}|^2,\\
			&=\int_{0}^{1}|e^{2\pi ij_1x}+e^{2\pi i j_2x}|^2dx,\\
			&=\left\{	\begin{aligned}
				2\quad\quad	&\text{if $j_1\neq j_2,$}\\
				4\quad\quad&\text{if $j_1 = j_2$.}
			\end{aligned}\right.
		\end{aligned}
	\end{equation*}
	Similarly,
	\begin{equation*}
		\begin{aligned}
			{Q}^2(T_{j_1k_1}-T_{j_2k_2})&=\left\{\begin{aligned}
				2\quad\quad\quad\quad\quad\quad\quad\quad	&\text{if $k_1\neq k_2,$}\\
				\lim_{n\to\infty}\frac{1}{n}\sum_{r=0}^{\lfloor \frac{n}{\lfloor \sqrt{n}\rfloor}\rfloor-1}\sum_{l\in I_{\lfloor \sqrt{n}\rfloor}}|e^{2\pi ij_1l}+e^{2\pi i j_2l}|^2\quad\quad&\text{if $k_1 = k_2$,}
			\end{aligned}\right.\\
			{Q}^2(T_{j_1k}-T_{j_2k})&=\left\{	\begin{aligned}
				2\quad\quad	&\text{if $j_1\neq j_2,$}\\
				0\quad\quad&\text{if $j_1 = j_2$.}
			\end{aligned}\right.
		\end{aligned}
	\end{equation*}
	Then,
	$$Re(\langle T_{j_1k_1}+Z, T_{j_2k_2}+Z\rangle)=\left\{	\begin{aligned}
		1\quad\quad	&\text{if $j_1= j_2, k_1=k_2$}\\
		0\quad\quad&\text{otherwise.}
	\end{aligned}\right.$$
	Also,
	$$Im(\langle T_{j_1k_1}+Z, T_{j_2k_2}+Z\rangle)=0.$$
	
\end{proof}
	The next two theorems provide the motivation for considering the seminorm $Q$ in defining the LT and GLT operator.
		\begin{theorem}\label{GLT}
		If $T$ is a GLT Operator, then $\{P_nTP_n\oplus O_{n \mod \lfloor \sqrt{n}\rfloor}\}_n$ is a GLT sequence.
	\end{theorem}
	\begin{proof}
		For $\epsilon>0$, there exist an integer $M$ such that, for all $m>M$,
		$$\|T- \sum_{j_m,k_m}a_{j_mk_m}T_{j_mk_m}\|_G<\epsilon.$$
		Let $P_{C_n}$ and $P_{D_n}$ are the projection from $\mathcal{H}^2$ onto the space spanned by $C_n$ and $D_n$ receptively, where  $C_n$ and $D_n$ are the decomposition of $B_n$. Then Rank$((T_n-T_{n,m})(P_{D_n}))=o(n)$. Now
		\begin{equation*}
			\begin{aligned}
				q_{2}^{2}(\{T_n-T_{n,m}\}_n)&\leq \limsup\frac{1}{n}\|(T_n-T_{n,m})P_{C_n}\|_F^2\\
				&\leq\limsup_{n\to\infty}\frac{1}{n}\sum_{e_{l,r}\in C_n}\|T- \sum_{j_m,k_m}a_{j_mk_m}T_{j_mk_m}\|^2\\
				&<\epsilon.
			\end{aligned}
		\end{equation*}	
	\end{proof}
\begin{theorem}\label{cvgn rln}
	Operator norm convergence $\implies$ $\|\cdot\|_G$ convergence $\implies$ $Q$ convergence.
\end{theorem}
\begin{proof}
	Straight forward from the definitions.
\end{proof}

	\section{LT Operator and GLT Operator}
	In this section, we develop the theory of LT Operators followed by GLT Operators.
	The results contained in this section are of fundamental importance for the theory of GLT Operators.	
		A class of linear operators that plays a central role here is the class of zero-distributed operator.
	
	\begin{definition}
		An operator  $Z:D(Z)\to \mathcal{H}^2$ is zero distributed operator if $Q(Z)=0.$
	\end{definition}
	Let us denote $Z(\mathcal{H}^2)$ be the set of all Zero distributed operators. Following are some examples for zero distributed operators.
	\begin{example}\label{Zerexa1}
		From the definition of Zero distributed operator it is clear that every compact operator on $\mathcal{H}^2$ is a zero distributed operator. 
	\end{example}
	\begin{example}\label{Zerpexa2}
		$T:\mathcal{H}^2\to\mathcal{H}^2$ be a bounded linear operator defined by $\displaystyle T(\sum_{l,k}x_{lk}e_{lk})=\sum_{l}x_{l1}e_{l1}$. Here $T$ is not a compact operator.
	\end{example}
	The examples \ref{Zerexa1} and \ref{Zerpexa2} says $K(\mathcal{H}^2)$, the space of compact operator on $\mathcal{H}^2$ is a proper subspace of $Z(\mathcal{H}^2)$.\\Now we give a characterization of  zero distributed operators.
	\begin{theorem}
		Let $Z\in L(\mathcal{H}^2)$, then $Z$ is zero distributed operator if and only if for every $\epsilon>0$, $\lim_{n\to\infty}\frac{1}{n}\#\{e_{l,r}\in B_n: \|Z(e_{l,r})\|>\epsilon\}=o.$
	\end{theorem}
	\begin{proof}
		Let $Z$ be a zero distributed operator. Suppose there exist  $\epsilon_0, \epsilon_1>0$ and a sub-sequence $\{n_k\}$ such that 
		$$\frac{1}{n_k}\#\{e_{l,r}\in B_{n_k}: \|T(e_{l,r})\|>\epsilon_0\}>{\epsilon_1}\quad \forall k.$$
		Let $C_0=\{e_{l,r}\in B_{n_k}: \|Z(e_{l,r})\|>\epsilon_0\}$. Then every decomposition $C_{n_k}\sqcup D_{n_k}=B_{n_k}$, $\lim_{n\to\infty}\frac{1}{n}\#(C_n\cap C_0)>\epsilon_1$. Then,
		$$\limsup_{n\to \infty}\frac{1}{n}\sum_{e_{l,r}\in C_n}\|Z(e_{l,r})\|^2>\epsilon_1\epsilon_0^2,$$
		which is a contradiction.\\
		The converse is trivial.
	\end{proof}
	The next theorem gives the relation between zero distributed operator and zero distributed  sequences.
	\begin{theorem}
		Let $Z$ be a zero distributed operator, then $\{P_nZP_n\oplus O_{n\, mod\lfloor\sqrt{n}\rfloor}\}$ is a zero distributed sequence. 
	\end{theorem}
	\begin{proof}
		Let $Z_n=P_nZP_n\oplus O_{n\, mod\lfloor\sqrt{n}\rfloor}$,
		\begin{equation*}
			\begin{aligned}
				\lim_{n\to\infty}\frac{1}{n}\|Z_n\|_F^2&\leq \lim_{n\to\infty}\frac{1}{n}\sum_{r=0}^{\lfloor \frac{n}{\lfloor \sqrt{n}\rfloor}\rfloor-1}\sum_{l\in I_{\lfloor \sqrt{n}\rfloor}}\|P_nZP_n(e_{lr})\|^2\\
				&\leq\lim_{n\to\infty}\frac{1}{n}\sum_{r=0}^{\lfloor \frac{n}{\lfloor \sqrt{n}\rfloor}\rfloor-1}\sum_{l\in I_{\lfloor \sqrt{n}\rfloor}}\|Z(e_{lr})\|^2\\
				&=0.
			\end{aligned}
		\end{equation*}
	By Lemma \ref{zerod} $\{Z_n\}_n$ is a zero distributed sequence.
\end{proof}
	If $ T$ be a bounded linear operator on $\mathcal{H}^2$ and $Z$ be a bounded zero distributed operator on $\mathcal{H}^2$, then
	\begin{equation}\label{ideal}
		\lim_{n\to\infty}\frac{1}{n}\sum_{r=0}^{\lfloor \frac{n}{\lfloor \sqrt{n}\rfloor}\rfloor-1}\sum_{l\in I_{\lfloor \sqrt{n}\rfloor}}\|TZ(e_{lr})\|^2dt\\
		\leq\lim_{n\to\infty}\frac{1}{n}\sum_{r=0}^{\lfloor \frac{n}{\lfloor \sqrt{n}\rfloor}\rfloor-1}\sum_{l\in I_{\lfloor \sqrt{n}\rfloor}}\|T\|^2\|Z(e_{lr})\|^2dt=0.
	\end{equation}	
	From (\ref{ideal}) we can see that the set of all bounded zero distributed operator is a left ideal in $B(\mathcal{H}^2)$.

Now we are ready to introduce the notion of LT and GLT Operators.
\begin{definition}
		An Operator $T:D(T)\to \mathcal{H}^2$ is said to be a LT operator if $T:=T_{af}+Z$, where $Z$ is a zero distributed operator.
\end{definition}
\begin{definition}
	An Operator $T:D(T)\to \mathcal{H}^2$ is said to be a GLT operator if there exist sequence of operators $\{\sum_{j_m,k_m}a_{j_mk_m}T_{j_mk_m}\}$ such that 
	$$\lim_{m\to\infty}Q(T- \sum_{j_m,k_m}a_{j_mk_m}T_{j_mk_m})=0.$$
\end{definition}
	From Theorem \ref{GLT} we understood that, if $T$ is a GLT operator such that $\sum_{j_m,k_m}a_{j_mk_m}T_{j_mk_m}$ converges to $T$, then there exist a measurable function $\kappa$ such that $\sum_{j_m,k_m}a_{j_mk_m}e^{2\pi j_mx}\cdot e^{ik_mt}$ converges to $\kappa$ in $L^p$ and $\{T_n\}_n$ is a GLT sequence with symbol $k$. In this case we say that $T$ is a GLT operator operator with symbol $k$ and it is denoted by $T\sim_{GLT}k$.\\
	We know from \cite{garoni2017generalized} every Locally Toeplitz sequence is a GLT sequence. The next theorem says the same is true for Locally Toeplitz Operator also, that is every Locally Toeplitz Operator is a GLT operator.
	\begin{theorem}
		$LT_{af}$ is a GLT operator
	\end{theorem}
	\begin{proof}
		Result follows from Theorem \ref{LT-GLT} and Theorem \ref{cvgn rln}
	\end{proof}
	\begin{corollary}
		Let $LT_{af}$ be a Locally Toeplitz Operator and $T_n=P_nLT_{af}P_n \oplus O_{n\mod \lfloor \sqrt{n}\rfloor}$. Then $\{T_n\}_n$ is a LT sequence. 
	\end{corollary}
	Since each $T_{jk}$ is LT operator and all LT operator is a GLT operator, the definition of GLT operator can be characterized in terms of LT operator just like as in the case of matrix sequences.
	\begin{theorem}
		An Operator $T:D(T)\to \mathcal{H}^2$ is said to be a GLT operator if there exist sequence of operators $\{\sum_{j,k}T_{a_jf_j}\}$ such that 
		$$\lim_{n\to\infty}\|T- \sum_{j,k}T_{a_jf_j}\|_G=0.$$	
	\end{theorem}
	Now we investigate the important algebraic properties of GLT operator.
	\begin{theorem}
		Suppose that $A\sim_{GLT}\kappa$ and $B\sim_{GLT}\zeta$. Then,
		\begin{enumerate}
			\item $\alpha A+\beta B\sim_{GLT}\alpha\kappa+\beta\zeta$.
		\end{enumerate}
	\end{theorem}

	\begin{remark}\normalfont\label{example}
		Product of two GLT operators may not be a GLT operator. For, take LT operator $LT_{af}$ with $a(x)=x^{\frac{-1}{4}}$ and $f=1$. Then as in the proof of Theorem \ref{onb}, we have to just calculate $\|\cdot\|_G$. 
		\begin{equation*}
			\begin{aligned}
				\|LT_{x^{\frac{-1}{4}}\cdot 1}\|_G^2&=\lim_{n\to\infty}\frac{1}{n}\sum_{r=0}^{\lfloor \frac{n}{\lfloor \sqrt{n}\rfloor}\rfloor-1}\sum_{l\in I_{\lfloor \sqrt{n}\rfloor}}\|LT_{x^{\frac{-1}{4}}\cdot 1}(e_{lr})\|^2\\
				&=\lim_{n\to\infty}\frac{1}{n}\sum_{r=0}^{\lfloor \frac{n}{\lfloor \sqrt{n}\rfloor}\rfloor-1}\sum_{l\in I_{\lfloor \sqrt{n}\rfloor}}\|l^{\frac{-1}{4}}e^{irt}\|^2\\	
				&=\lim_{n\to\infty}\frac{1}{n}\sum_{r=0}^{\lfloor \frac{n}{\lfloor \sqrt{n}\rfloor}\rfloor-1}\sum_{l\in I_{\lfloor \sqrt{n}\rfloor}}\frac{1}{2\pi}\int_{-\pi}^{\pi}l^{\frac{-1}{2}}dt
			\end{aligned}
		\end{equation*}
	\begin{equation*}
\begin{aligned}
				&=\lim_{n\to\infty}\frac{1}{\lfloor \sqrt{n}\rfloor}\sum_{l\in I_{\lfloor \sqrt{n}\rfloor}}l^{\frac{-1}{2}}\\
				&=\int_0^1 x^{\frac{-1}{2}}dx\\
				&=2	
			\end{aligned}
		\end{equation*}
	\begin{equation*}
		\begin{aligned}
				\|LT_{x^{\frac{-1}{4}}\cdot 1}^2\|_G^2&=\lim_{n\to\infty}\frac{1}{n}\sum_{r=0}^{\lfloor \frac{n}{\lfloor \sqrt{n}\rfloor}\rfloor-1}\sum_{l\in I_{\lfloor \sqrt{n}\rfloor}}\|LT_{x^{\frac{-1}{4}}\cdot 1}^2(e_{lr})\|^2\\
				&=\lim_{n\to\infty}\frac{1}{n}\sum_{r=0}^{\lfloor \frac{n}{\lfloor \sqrt{n}\rfloor}\rfloor-1}\sum_{l\in I_{\lfloor \sqrt{n}\rfloor}}\|LT_{x^{\frac{-1}{4}}\cdot 1}(l^{\frac{-1}{4}}e_{lr})\|^2\\
				&=\lim_{n\to\infty}\frac{1}{n}\sum_{r=0}^{\lfloor \frac{n}{\lfloor \sqrt{n}\rfloor}\rfloor-1}\sum_{l\in I_{\lfloor \sqrt{n}\rfloor}}\|l^{\frac{-1}{2}}e^{irt}\|^2\\
				&=\lim_{n\to\infty}\frac{1}{n}\sum_{r=0}^{\lfloor \frac{n}{\lfloor \sqrt{n}\rfloor}\rfloor-1}\sum_{l\in I_{\lfloor \sqrt{n}\rfloor}}\frac{1}{2\pi}\int_{-\pi}^{\pi}l^{-1}dt\\
				&=\lim_{n\to\infty}\frac{1}{\lfloor \sqrt{n}\rfloor}\sum_{l\in I_{\lfloor \sqrt{n}\rfloor}}l^{-1}\\
				&=\int_0^1 x^{-1}dx\\
				&=\infty
			\end{aligned}
		\end{equation*}
		Thus $LT_{x^{\frac{-1}{4}}\cdot 1}$ is a GLT operator but the product $LT_{x^{\frac{-1}{4}}\cdot 1}^2$ is not a GLT.
	\end{remark}

	 The next next theorem is analogous to theorem in \cite{sarath}, says that the space of all GLT operator is isomorphic to $L^2([0,1]\times[-\pi,\pi])$ as Banach space.
	
	\begin{theorem}	$\mathcal{G}$ and $L^2([0,1]\times[-\pi,\pi])$ are isometrically isomorphic. 
	\end{theorem}
	\begin{proof}
		Let $T\in E$, that is $T=\sum_{j,k}a_{jk}T_{jk}$ and $F$ is the span of $\{e^{2\pi x}e^{ikt}:j,k\in \mathbb{Z}\}$.
		\begin{equation*}
			\begin{aligned}
				\|T\|_G^2&=\lim_{n\to\infty}\frac{1}{n}\sum_{r=0}^{\lfloor \frac{n}{\lfloor \sqrt{n}\rfloor}\rfloor-1}\sum_{l\in I_{\lfloor \sqrt{n}\rfloor}}\frac{1}{2\pi}\int_{-\pi}^{\pi}|\sum_{j,k}a_{jk}T_{jk}(e_{lr})|^2 dt,\\
				&=\lim_{n\to\infty}\frac{1}{n}\sum_{r=0}^{\lfloor \frac{n}{\lfloor \sqrt{n}\rfloor}\rfloor-1}\sum_{l\in I_{\lfloor \sqrt{n}\rfloor}}\frac{1}{2\pi}\int_{-\pi}^{\pi}|\sum_{j,k}a_{jk}e^{2\pi ijl}\cdot e^{i(k+r)t}|^2 dt,\\
				&=\lim_{n\to\infty}\frac{1}{n}\sum_{r=0}^{\lfloor \frac{n}{\lfloor \sqrt{n}\rfloor}\rfloor-1}\sum_{l\in I_{\lfloor \sqrt{n}\rfloor}}\frac{1}{2\pi}\sum_{k}|\sum_{j}a_{jk}e^{2\pi ijl}|^2,
			\end{aligned}
		\end{equation*}
	\begin{equation*}
\begin{aligned}
				&=\sum_{k}\lim_{n\to\infty}\frac{1}{\lfloor \sqrt{n}\rfloor}\sum_{l\in I_{\lfloor \sqrt{n}\rfloor}}|\sum_{j}a_{jk}e^{2\pi ijl}|^2,\\
				&=\sum_{k}\int_{0}^{1}|\sum_{j}a_{jk}e^{2\pi ijx}|^2 dx,\\
				&=\sum_{j,k}|a_{jk}|^2.
			\end{aligned}
		\end{equation*}
		Now the map $\displaystyle \phi:E\to F$, $\phi(\sum_{j,k}a_{jk}T_{jk})=\sum_{j,k}a_{jk}e^{2\pi ijx}\cdot e^{ikt}$ is an isometric isomorphism of $E$ onto $F$. Then by the continuity of norm $\mathcal{G}$ and $L^2([0,1]\times[-\pi,\pi])$ are isometrically isomorphic.
	\end{proof}
\begin{corollary}\label{symbol}\sloppy
	Let $T\in \mathcal{G}$ with symbol {$f\in L^2([0,1]^d\times[-\pi,\pi]^d)$}, let $\displaystyle \sum_{\textbf{r},\textbf{s}=-\infty}^{\infty} a_{\textbf{jk}}e^{i2\pi \textbf{j}\cdot{x}}e^{i\textbf{k}\cdot {t}}$  be the  Fourier series representation of $f$. Then	$	a_{\textbf{jk}}= \langle T, T_{jk}\rangle_{\mathcal{G}},$ 
	where $\langle \cdot, \cdot \rangle_{\mathcal{G}}$ is the inner product induced from $Q(T)$.

\end{corollary}
\begin{remark}
	If we choose matrix sequence $\{T_n\}_n$ instead of $T$, $\lim_{n\to \infty}\langle T_n,T_{jk}\rangle$ coincide with the equation defined in Theorem 2.2 of \cite{Sarath2024}.
\end{remark}
The following result is important in the application described in the next section.
\begin{theorem}\label{comsymbol}
		If both T and S belong to $\mathcal{G}$ and {suppose that}  for all $D_n\subset B_n$ such that $|D_n|=o(n)$ and $\lim_{n\to \infty}\frac{1}{n}\sum_{e_{lk}\in D_n}\|T(e_{lk})\|^2=\lim_{n\to \infty}\frac{1}{n}\sum_{e_{lk}\in D_n}\|S(e_{lk})\|^2=0$, then
		$$\langle T,S\rangle_\mathcal{G}= \lim_{n\to \infty}\frac{1}{n}\sum_{e_{lr}\in B_n}\langle T(e_{lk}),S(e_{lr})\rangle.$$
\end{theorem}

	\begin{proof}
		From the hypothesis we obtain, $Q(T)=\limsup_{n\to \infty}\frac{1}{\sqrt{n}}\left(\sum_{e_{lk}\in B_n}\|T(e_{lk})\|^2\right)^{\frac{1}{2}}$.
		Since $T\in \mathcal{G}$,  there exist sequence $\{T_m\}$ of $E$ such that $Q(T-T_m)<\epsilon.$ Now,
		\begin{equation*}
			\begin{aligned}
				Q(T)&<Q(T_m)+Q(T-T_m),\\
				\limsup_{n\to \infty}\frac{1}{\sqrt{n}}\left(\sum_{e_{lk}\in B_n}\|T(e_{lk})\|^2\right)^{\frac{1}{2}}&\leq \limsup_{n\to \infty}\frac{1}{\sqrt{n}}\left(\sum_{e_{lk}\in B_n}\|T_m(e_{lk})\|^2\right)^{\frac{1}{2}}+Q(T-T_m)\\,
				&\leq \liminf_{n\to \infty}\frac{1}{\sqrt{n}}\left(\sum_{e_{lk}\in B_n}\|T_m(e_{lk})\|^2\right)^{\frac{1}{2}}+Q(T-T_m),\\
				&\lim_{m\to\infty}\liminf_{n\to \infty}\frac{1}{\sqrt{n}}\left(\sum_{e_{lk}\in B_n}\|T_m(e_{lk})\|^2\right)^{\frac{1}{2}}+\lim_{m\to\infty}Q(T-T_m),\\
				&\liminf_{n\to \infty}\frac{1}{\sqrt{n}}\left(\sum_{e_{lk}\in B_n}\|T(e_{lk})\|^2\right)^{\frac{1}{2}}.
			\end{aligned}
		\end{equation*}
	Thus $\displaystyle Q(T)=\lim_{n\to \infty}\frac{1}{\sqrt{n}}\left(\sum_{e_{lk}\in B_n}\|T(e_{lk})\|^2\right)^{\frac{1}{2}}.$ Now from polarization identity, 
	$$\langle T,S\rangle_\mathcal{G}= \lim_{n\to \infty}\frac{1}{n}\sum_{e_{lk}\in B_n}\langle T(e_{lk}),S(e_{lk})\rangle.$$	
	\end{proof}	
	\section{Application}
	Let $T$ is a GLT operator, then $\{T_n\}_n$ is a GLT sequence and $T(e_{l,r})$ is approximately equal to $ T_n(e_{lr})$. Then Theorem \ref{comsymbol} and Corollary \ref{symbol} enable us to find the symbol of some GLT sequences.
	Consider the second-order differential problem:
	\begin{equation}
		 \begin{aligned}\label{eqn2}
				-(a(x)u^\prime(x))^\prime=f(x),\quad &x\in (0,1),\
				u(0)=\alpha,\quad u(1)=\beta,
			\end{aligned}
		\end{equation}
		where $a:[0,1]\to \mathbb{R}$ is continuous. we explore the discretization approach employing the classical second-order central FD scheme. We opt for a discretization parameter $n\in \mathbb{N}$ and set $x_j=jh$ for all $j\in [0,n+1]$, where $h=\frac{1}{n+1}$. It's noteworthy that, for $j=1,2,\ldots,n$, we can estimate $ -({a(x)}u^\prime(x))^\prime|_{x=x_j}$ using the subsequent FD formula:
	
	\begin{equation}\label{eqn3}
		\begin{aligned}
			-(a(x)u^\prime(x))^\prime&\approx-\frac{a(x_{j+\frac{1}{2}})u^\prime (x_{j+\frac{1}{2}})-a(x_{j-\frac{1}{2}})u^\prime(x_{j-\frac{1}{2}})}{h}\\&\approx-\frac{1}{h}\left(a(x_{j+\frac{1}{2}})\frac{u(x_{j+1})-u(x_j)}{h}-a(x_{j-\frac{1}{2}})\frac{u(x_j)-u(x_{j-1})}{h}\right)\\
			&=\frac{1}{h^2}\left(-a(x_{j+\frac{1}{2}})u(x_{j+1})+(a(x_{j+\frac{1}{2}})+a(x_{j-\frac{1}{2}}))u(x_j)-a(x_{j-\frac{1}{2}})u(x_{j-1})\right).
		\end{aligned}
	\end{equation}
	Subsequently, we approximate the solution of (\ref{eqn2}) using a piecewise linear function that is globally continuous, taking the value $u_j$ at $x_j$ for $j=0,\ldots,n+1$, where $u_0=\alpha$, $u_{n+1}=\beta$, and $u=(u_1,\ldots,u_n)^T$ denotes the solution of the linear system
	\begin{equation}\label{eqn4}
		-a(x_{j+\frac{1}{2}})u_{j+1}+(a(x_{j+\frac{1}{2}})+a(x_{j-\frac{1}{2}}))u_j-a(x_{j-\frac{1}{2}})u_{j-1}=h^2f(x_j),\quad j=1,\ldots,n.
	\end{equation}
	The matrix representing the linear system (\ref{eqn4}) is the tridiagonal symmetric matrix:
	\begin{equation*}
		A_n=\begin{pmatrix}
			a(x_{\frac{1}{2}})+a(x_{\frac{3}{2}})&-a(x_{\frac{3}{2}})&&&&\\
			-a(x_{\frac{3}{2}})&a(x_{\frac{3}{2}})+a(x_{\frac{5}{2}})&-a(x_{\frac{5}{2}})&&&\\
			&-a(x_{\frac{5}{2}})&a(x_{\frac{5}{2}})+a(x_{\frac{7}{2}})&-a(x_{\frac{7}{2}})&&\\
			&&-a(x_{\frac{7}{2}})&\ddots&\ddots&\\
			&&&\ddots&\ddots&-a(x_{n-\frac{1}{2}})\\
			&&&&-a(x_{n-\frac{1}{2}})&a(x_{n-\frac{1}{2}})+a(n+x_{\frac{1}{2}})
		\end{pmatrix}.
	\end{equation*}
This example demonstrates that Theorem \ref{comsymbol} enables the computation of the singular value and eigenvalue distribution of the sequence of discretization matrices ${A_n}$. Utilizing (\ref{innerprrdct}), we can evaluate the function $\displaystyle f_l=\sum_{r,s=-l}^{l} a_{rs}e^{i2\pi rx}e^{is\theta}$. Consequently, $\displaystyle\lim_{l\to\infty}f_l=f$ represents the symbol of $\{A_n\}_n$. Hence, $\{A_n\}_n\sim_{GLT}f$, and given the Hermitian nature of $A_n$, $f$ must be a real-valued function, establishing $\{A_n\}_n\sim_{\lambda}f$.\\
Consider $a(x)=2 \sin(x)+\cos(2x)$.  A simple computation shows that {$\displaystyle a_{rs}=\lim_{n\to \infty}\frac{1}{n}\sum_{e_{p,q}\in B_n}\langle T_m(e_{p,q}),S(e_{p,q})\rangle=0$, for every $s\in \mathbb{Z}\setminus\{-1,0,1\}$.\\
 Now, for $l=3,n=40^2$ and $m=1000$,
\begin{equation}\label{symbolf}
	f_3=
	\begin{bmatrix}
		e^{-6i\pi x}\\
		e^{-4i\pi x}\\
		e^{-2i\pi x}\\
		1\\
		e^{2i\pi x}\\
		e^{4i\pi x}\\
		e^{6i\pi x}\\
	\end{bmatrix}^T
	\begin{bmatrix}
		0.0089+i0.0135& -0.0153-i0.0271&0.0087+i0.0137\\
		0.0188+i0.0191&-0.0351-i0.0383&0.0186+i0.0193\\
	    0.0763+i0.0238& -0.1500-i0.0480& 0.0762+i0.0243\\
		-1.3744&2.7484&-1.3744\\
		0.0763-i0.0238&-0.1500+i0.0480&0.0762-i0.0243\\
		0.01881-i0.0191&-0.0351+i0.0383&0.0186-i0.0193\\
	    0.0089-i0.0135&-0.0153+i0.0271&0.0087-i0.0137\\
	\end{bmatrix}
	\begin{bmatrix}
		e^{-i\theta}\\
		1\\
		e^{i\theta}\\
	\end{bmatrix}.
\end{equation}
For $l=3$ and $n=400$, Table 1 displays the computation of the norm of $\gamma_n$ for increasing values of $m$. Here, $\gamma_n$ represents the vector comprising samples $Im(f_l(\frac{j}{\sqrt{n}},\frac{2\pi k}{\sqrt{n}}))$, where $j,k=0,1,\ldots,\sqrt{n}-1$. Given the negligible nature of $Im(f_l)$, focusing solely on the real part of $f_l$ suffices.
\begin{table}[ht]
	\centering	
	\begin{tabular}{|c|c c c c|}
		\hline
		{$m$} & 100 & 400 & 700& 1000 \\
		\hline
		$\|\gamma_n\|_2$&0.01104&0.00214&0.00119&0.00822\\
		\hline
	\end{tabular}
	\caption{Computation of $\|\gamma_n\|_2$ for increasing values of $m$.}
\end{table}
Table 2 presents the 2-norm of the difference $\zeta_n-\eta_n$ for increasing values of $n$ and $l$, where:
\begin{itemize}
	\item $\zeta_n$ denotes the vector of eigenvalues of $A_n$.
	\item $\eta_n$ represents the vector of samples $Re({f_l}(\frac{j}{\sqrt{n}},\frac{2\pi k}{\sqrt{n}}))$, with $j,k=0,1,\ldots,\sqrt{n}-1$.
\end{itemize}
Both $\zeta_n$ and $\eta_n$ are sorted in non-increasing order.

From \cite{garoni2017generalized}, we know that the symbol of ${A_n}_n$ is $f(x,\theta)=a(x)(2-2\cos(\theta))$. Table 3 presents the computation of the 2-norm of the function $f-f_l$ for increasing values of $m$ and $l$.

In Figure 2, we depict the spectrum of $A_n$ alongside the values $Re({f_l}(\frac{j}{\sqrt{n}},\frac{2\pi k}{\sqrt{n}}))$ and $f(\frac{j}{\sqrt{n}},\frac{2\pi k}{\sqrt{n}})$ for $n=3600$ and $l=10$. The eigenvalues of $A_n$ and the sampling of $Re(f_l)$ and $f$ are shown in non-increasing order.
\begin{figure}[h!]
	\centering
	\includegraphics{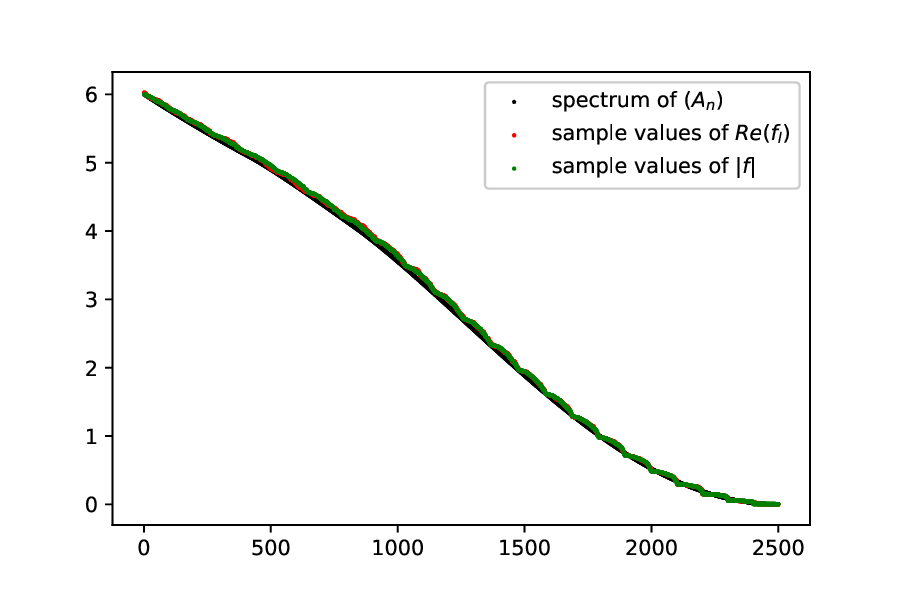}
	\caption{Spectrum of $A_n$ and samplings $Re({f_l}(\frac{j}{\sqrt{n}},\frac{2\pi k}{\sqrt{n}})),f(\frac{j}{\sqrt{n}},\frac{2\pi k}{\sqrt{n}}),j,k=0,1,\ldots,\sqrt{n}-1,$ for $n=50^2$ and $l=7$.}
\end{figure}
\begin{table}[h!]
	\centering	
	\begin{tabular}{|c|c c c c |}
		\hline
		\backslashbox{$l$}{$n$} & $20^2$ & $40^2$ & $60^2$ & $80^2$ \\
		\hline
		5&0.1559&0.0735&0.0489&0.0364\\
		\hline
		10&0.1524&.0712&0.0466&0.0347\\
		\hline
		15&0.1522&0.0706&0.0461&0.0343\\
		\hline
	\end{tabular}
	\caption{Computation of $\|\zeta_n-\eta_n\|_2$ for increasing values of $n$ and $l$.}
\end{table}
\begin{table}[h!]
	\centering	
	{\begin{tabular}{|c|c c c c|}
			\hline
			\backslashbox{$l$}{$m$} & $100$ & $400$ & $700$  & $1000$\\
			\hline
			3&0.0889&0.0809&0.0805&0.0804\\
			\hline
			7&0.0784&0.0558&0.0545&0.0543\\
			\hline
			10&0.0816&0.0484&0.0465&0.0459\\
			\hline
	\end{tabular}}
	\caption{{Computation of $\|f-f_l\|_2$ for increasing values of $m$ and $l$.}}
\end{table}
\section{Conclusions}
This article's main accomplishment is that we tried to offer a foundation for creating an operator-theoretical version of GLT sequences. As an application of our main results, we proposed a technique to find the spectral symbol of GLT sequences under some mild conditions. The operator developed here has some drawbacks, such as the fact that the space of all GLT operators is not closed under multiplication and is not *-closed. Future research could thus focus on creating a new operator that gets around these problems.
\section*{Acknowledgment}
V. B. Kiran Kumar is thankful to KSCSTE, Kerala for financial support through the KSYSA-Research Grant. N.S. Sarathkumar is thankful to CSIR for financial support through CSIR-SRF.

\end{document}